\documentclass[12pt,a4paper]{article}
\usepackage[dvips]{graphicx}
\usepackage{setspace}
\usepackage{amsmath}
\usepackage{overcite}
\usepackage{amsbsy}
\usepackage{amssymb}
\usepackage[utf8]{inputenc}
\usepackage[T1]{fontenc}

\newcommand{\bmi}{\ensuremath{\boldsymbol}}

\newcommand{\bnabla}{\ensuremath{\boldsymbol\nabla}}
\newcommand{\bcdot}{\ensuremath{\boldsymbol\cdot}}

\newcommand{\beq}{\begin{equation}}
\newcommand{\eeq}{\end{equation}}

\long\def\symbolfootnote[#1]#2{\begingroup%
\def\thefootnote{\fnsymbol{footnote}}\footnote[#1]{#2}\endgroup}

\begin{document}

\begin{center}
{\bf \Large Two dimensional penetrative phototactic bioconvection with periodic sidewalls}\\[12pt]
{M. K. Panda$^1$\symbolfootnote[1]{Corresponding author;
e-mail:mkpanda@iiitdmj.ac.in}}
\\{$^1$ Department of Mathematics, PDPM Indian Institute of Information Technology Design and Manufacturing, Jabalpur 482005, India}
\end{center}

\noindent
{\bf Abstract}

Light gradient can allow many motile photosynthetic microorganisms to bias their motion towards moderate light (positive phototaxis) or away from intense light (negative phototaxis). The proposed work presents the penetrative phototactic bioconvection in a non-scattering algal suspension. The suspension is confined by a stress-free top boundary, and rigid bottom and periodic lateral boundaries. The resulting bioconvective patterns of the problem  strongly resemble to that of a spatially extended domain in the same vicinity. The bioconvection solution appears in the form of a two-rolls pattern (or any even number of rolls) due to the periodic lateral boundaries.

\newpage

\section{INTRODUCTION}
\label{sec1}

Bioconvection is recognized as the convective flows generated by the upward swimming motile microorganisms and it arises due to the differences in density between them and the fluid in their local environment \cite{pk:kp,hp:ph}. Usually, it appears spontaneously in biological suspensions and identified by emergent flow patterns.   When the microorganisms are of non-motile type, the flow patterns also do not appear. It is not difficult to find out examples of pattern formation where the microorganisms do not swim upward and are not heavier than the surrounding fluid \cite{pk:kp}. The swimming microorganisms found to bioconvect in biological suspensions are mostly bacteria, algae, and protozoa and so on. The response of microorganisms to external gradients or passive forces, which allow them to orient in particular directions in their environment is called \textit{taxes}. For instance, \textit{gravitaxis} represents a motion induced along the direction of gravity for a microorganism with asymmetry mass distribution.  In addition, if orientation of a microorganism is affected by a local shear due to fluid flow then the induced motion is named as \textit{gyrotaxis}. \textit{Phototaxis} is a mechanism by which a microorganism detects light gradients. This article accounts phototaxis only.

Laboratory experiments have shown that illumination can affect the bioconvection patterns \cite{wa,ka85,vin95}. For instance, the patterns in biological suspensions may be destroyed or their formation may be blocked via intense light. The reasons behind the changes in bioconvective flow patterns via light intensity may be as follows: The motile phototactic microorganisms detect light gradients and they exhibit positive phototaxis (i.e. $I<I_c$) and negative phototaxis (i.e. $I>I_c$). Thus, they select an optimal location across the suspension where $I\approx I_c$. The second reason is as follows. The algae absorb the light incident on them and it decreases along the incident direction of light. If the algal suspension is dilute and non-scattering, then the Lambert--Beer law gives the governing equation for light intensity $I(x,y,z)$ as \cite {vh:hv} i.e.

\begin{equation}
I(x,y,z)=I_t\, \exp{\left(\alpha \int_{\gamma}^{}n \,ds\right)}.
\label{intensity1}
\end{equation}

Here $I_t$ is the magnitude of vertical collimated solar flux, $\alpha$ is the absorption coefficient, $n$ is the cell concentration and 
$\gamma$ is a straight line segment connecting algae and the light source. If $<{\bmi p}>$ denotes the average swimming direction of cells in a small volume, then \cite{gh}

\begin{equation}
<{\bmi p}>=M(I) \hat{\bmi z},
\end{equation}
where  the swimming direction of a cell is defined by ${\bmi p},$ $\hat{\bmi z}$ is the unit vector directed vertically upwards. $M(I)$ is the photoresponse curve such that

\begin{equation}
M(I) \left\{
\begin{array} {c }
\ge 0 \quad \mbox{if} \quad I \le I_c, \\
< 0 \quad \mbox{if}\quad I > I_c.
\end{array}
\right.
\end{equation}

Here ${\bmi W_c} = W_c <{\bmi p}>$ represents the mean swimming velocity, where $W_c$ is the ensemble-average swimming speed.

Consider a biological suspension with an vertical collimated solar flux as the illuminating source. The equilibrium state for such a finte-depth suspension is the resultant between phototaxis with cell diffusion. As a result, a sublayer of algae is formed across the suspension as a function of intensity of light. The unstable zone is defined as the region below the sublayer and the region above it is defined as gravitationally unstable zone. The motion induced by convection from the unstable zone penetrates the upper stable zone \cite{bst}.

Many motile algae are strongly phototactic as they need to photosynthesize in order to survive. Thus, realistic and reliable models of their behaviour should include phototaxis. Phototaxis and the corresponding bioconvection also have important applications in biofuel, biohydrogen production, and biofluid dynamics research \cite{pg,Ishika-2009,williams_11}.

Gyrotactic, gravitactic and chemotactic bioconvection have been explored extensively via numerical simulation \cite{hp:ph,Lee_2015}. However, there has not been much numerical simulation on phototactic bioconvection available till date except few \cite{pg,gh,pr}. Two-dimensional phototactic bioconvection in a layer confined by a rigid bottom boundary, stress-free top and lateral boundaries was investigated by Ghorai and Hill \cite{gh} in the nonlinear regime. They used the phototaxis model proposed by Vincent and Hill \cite{vh:hv}. Panda and Ghorai \cite{pg} simulated numerically phototactic bioconvection in an isotropic scattering  suspension using the phototaxis model proposed by Ghorai \textit{et al.} \cite{gph}. Two-dimensional phototactic bioconvection was simulated numerically by Panda and Singh \cite{pr} in a layer confined by a stress-free top boundary, rigid bottom and lateral boundaries. Bioconvection in  small  domains  with  periodic  sidewalls has been hitherto carried out as a model of a spatially extended system and the proposed work investigates penetrative phototactic bioconvection  on the same vicinity. 

The manuscript is organized as follows: In Sec.~\ref{chapt4:mathformulation}, the problem is converted into mathematical form. In 
Sec.~\ref{chapt4:numerical_procedure}, the numerical procedure is proposed. In Sec.~\ref{chapt4:results}, the numerical results are presented. Finally, conclusions are drawn in Section~\ref{chapt4:conclusions}.

\section{Mathematical formulation}
\label{chapt4:mathformulation}

\subsection{Geometry of the problem}
\label{chapt4:geometry}

\begin{figure}[!h]
\begin{center}
\includegraphics[scale=0.65]{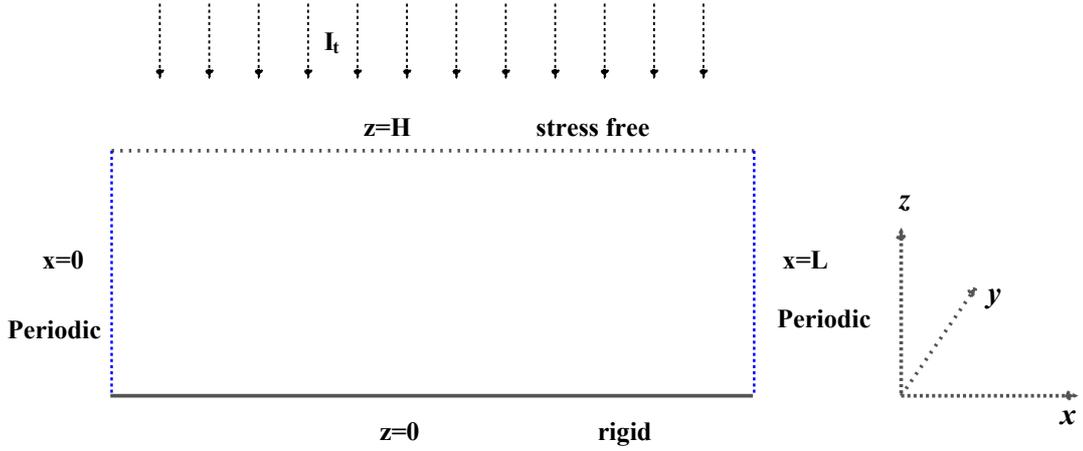}
\end{center}
\caption{Geometry of the problem. \label{chapt4:cdomain}}
\end{figure}

The geometry consists of a chamber with width $L$ and height $H$ referred to
Cartesian coordinates $\left(x,z\right)$ with the $z$-axis pointing vertically
upwards. Thus, the flow is confined in the $xz$-plane and is independent of 
the $y$ coordinate [see Fig.~\ref{chapt4:cdomain}]. 

\subsection{Governing equations}
\label{chapt4:governing}

We employ a continuum approach similar to the previous models on bioconvection \cite{vh:hv,pk:kp} and each cell has a volume $\vartheta$ and density $\rho+\Delta\rho$, where $\Delta\rho\ll \rho$ ad $\rho$ is the constant density of the water. Let ${\bmi u}=u\hat{\bmi x}+w\hat{\bmi z}$ and $n$ denote the average velocity and cell concentration over an elemental volume, where $\hat{\bmi x}$ and $\hat{\bmi z}$ are unit vectors along $x$ and $z$ axes. Assume that the suspension is incompressible and the flow is restricted to two dimensions only.  Introducing stream function $\psi$ and vorticity $\zeta$, we get

\begin{equation}
{\bmi u}=\left(\frac{\partial{\psi}}{\partial{z}},0,-\frac{\partial{\psi}}{\partial{x}}\right), \quad \zeta=-\nabla^2 \psi.
\end{equation}
 
The momentum equation under the Boussinesq approximation leads to the vorticity equation

\begin{equation}
\frac{\partial{\zeta}}{\partial{t}}+\bnabla \bcdot \left(\zeta {\bmi u}\right)=\nu \nabla^2 \zeta-\frac{\Delta \rho g \vartheta}{\rho}\,\frac{\partial{n}}{\partial{x}}.
\end{equation}

The conservation equation for microorganisms is given by

\begin{equation}
\frac{\partial{n}}{\partial{t}}=-{\bmi \nabla}\cdot{\bmi J},
\label{chapt4:scaled_cellflux0}
\end{equation}

where the flux of the cells is

\begin{equation}
{\bmi J}=n{\bmi u}+nW_{c}<{\bmi p}>-D {\bmi \nabla}{n}.
\label{chapt4:nondimensional_flux0}
\end{equation}

Here the first term,  second term and  third term represent the flux induced by advection of bulk fluid flow, average swimming of cells and random motion of cells respectively. We assume ${\textrm{\bf D}} = D{\textrm{\bf I}}.$  For a uniformly illuminated suspension via a vertical collimated solar flux as considered by Ghorai and Hill~\cite{gh}, the Eq.~(\ref{intensity1}) for light intensity becomes 

\begin{equation}
I(x,z)=I_t\,\exp{\left[-\alpha\,\int_{z}^{H} n(x,z) dz\right]}.
\end{equation}

\subsection{Boundary conditions}
\label{chapt4:boundary}

We impose rigid, no-slip boundary condition on $z=0$ and periodic boundary condition at $x=0, L$. We require that both the normal velocity and tangential stress vanish on $z=H$. Also there is no flux of cells through the walls. Thus the boundary conditions are

\begin{eqnarray}
 \psi=0 \quad \mbox{and\ }\quad {\bmi J}\cdot{\hat{\bmi z}}=0 \quad \mbox{at\ }\quad z=0,H,\\
 \frac{\partial{\psi}}{\partial{z}}=0 \quad \mbox{at}\quad z=0\quad \mbox{and}\quad \zeta=0\quad \mbox{at}\quad z=H,\\
\psi=0,~ \zeta=0 \quad \mbox{and}\quad {\bmi J}\cdot{\hat{\bmi x}=0 \quad \mbox{at}\quad} x=0,L,
\end{eqnarray}
and a periodic boundary condition with period $L$ is imposed in the $x-$ direction (i.e. $\psi{(0,z)}=\psi{(L,z)}$).

\begin{figure}[!h]
\begin{center}
\includegraphics[scale=1]{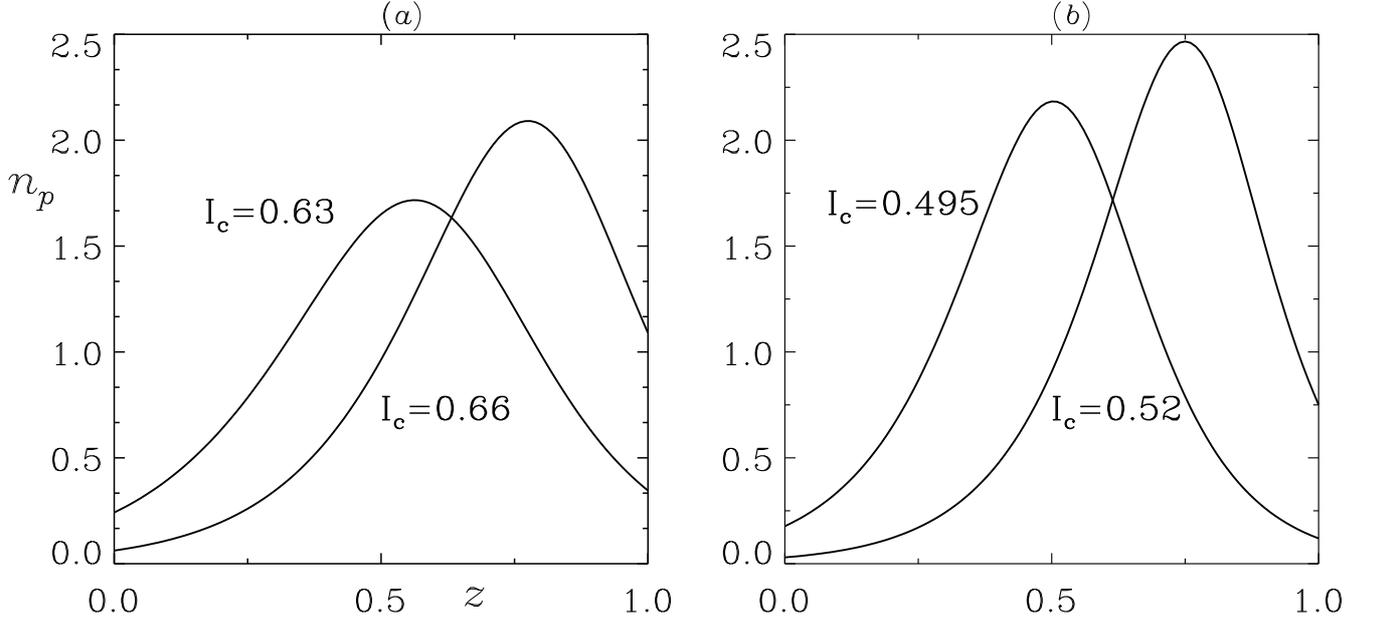}
\end{center}
\caption{Basic Equilibrium State   \label{basic}}
\end{figure}

\begin{figure}[!h]
\begin{center}
\includegraphics[scale=1]{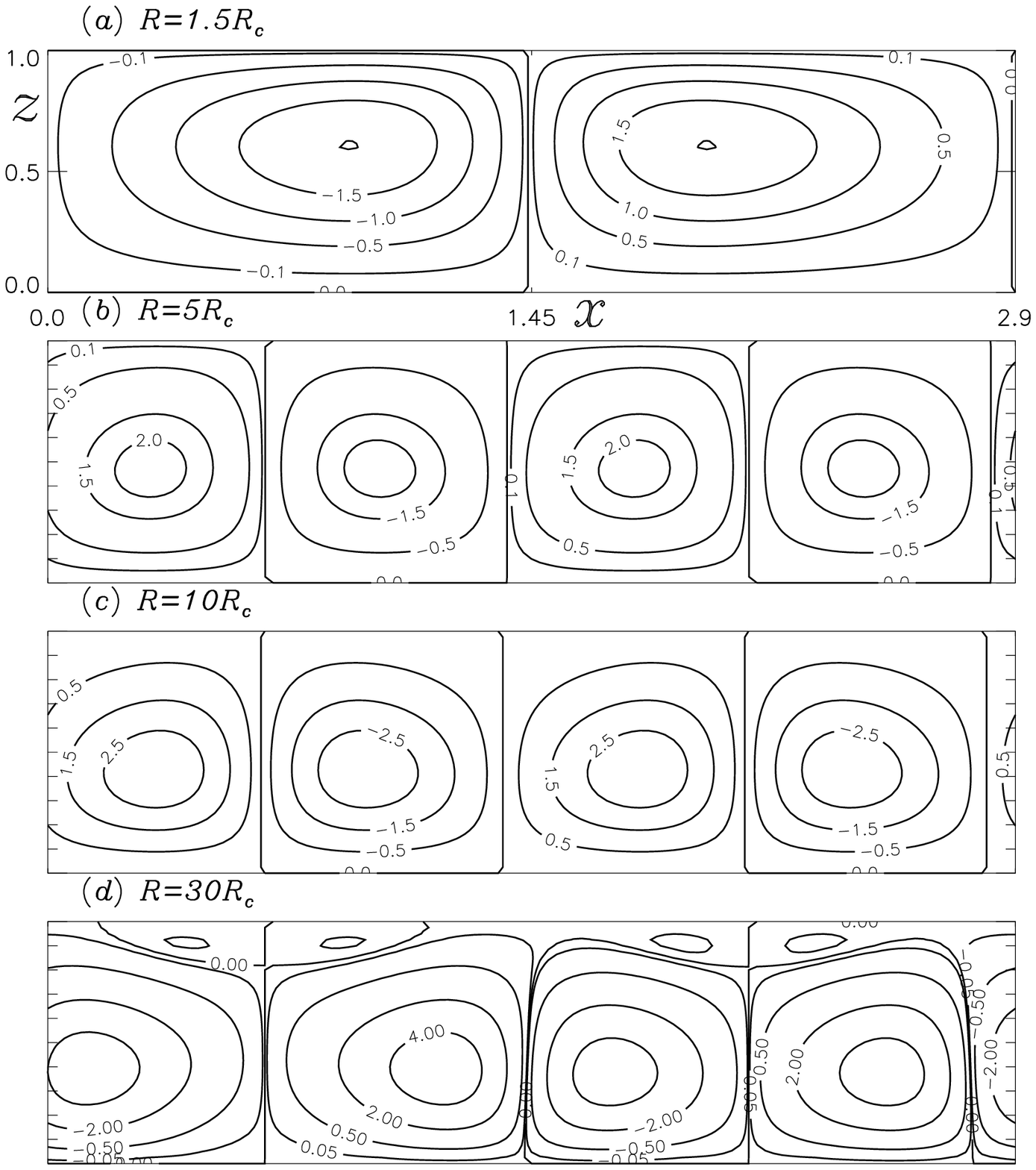}
\end{center}
\caption{Streamlines of steady solutions for $V_c=10$,~$\kappa=0.5$,~$I_c=0.66$. \label{Fig1}}
\end{figure}

\begin{figure}[!h]
\begin{center}
\includegraphics[scale=1]{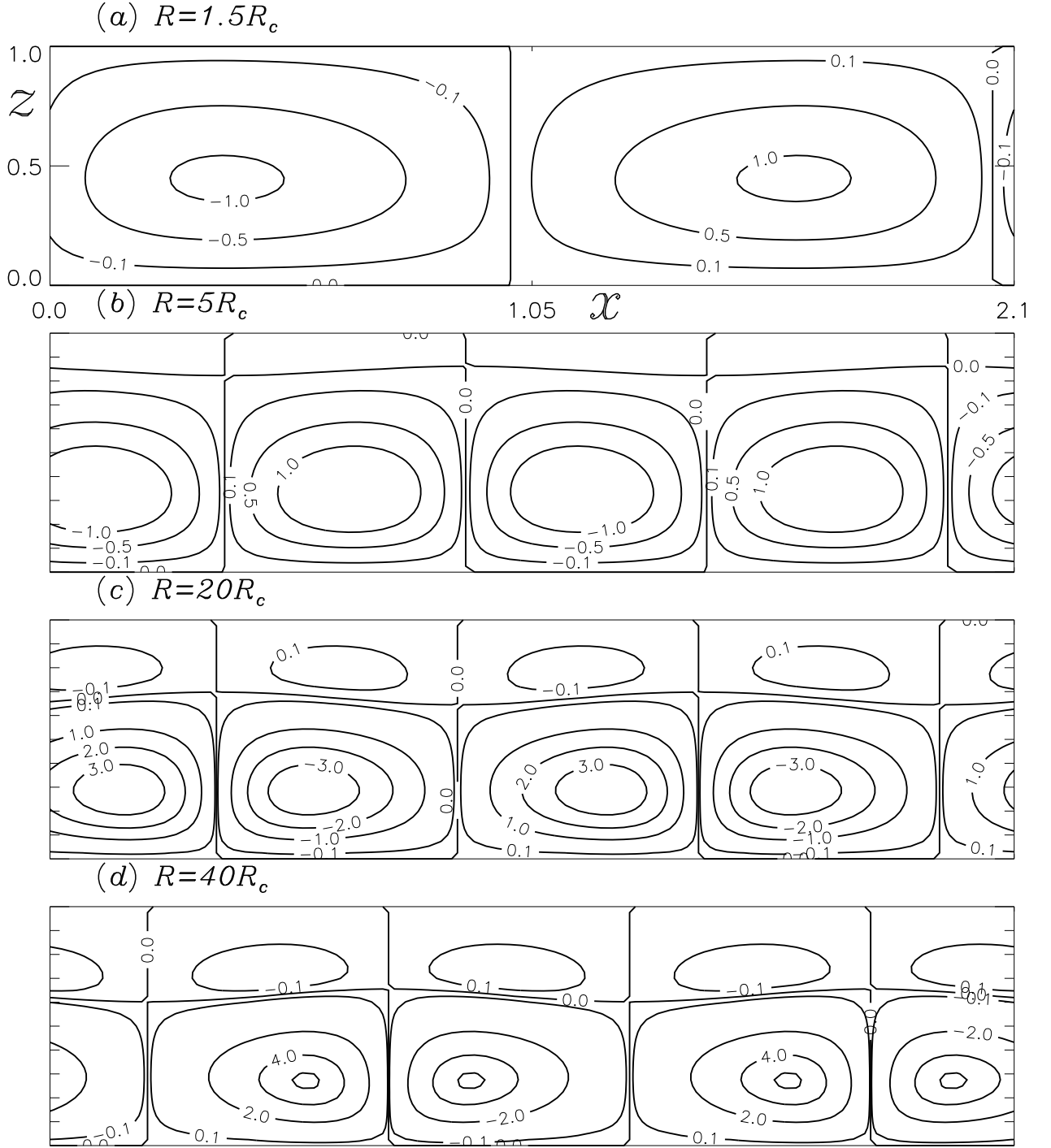}
\end{center}
\caption{Streamlines of steady solutions for $V_c=10$,~$\kappa=0.5$,~$I_c=0.63$. \label{Fig2}}
\end{figure}
 
\begin{figure}[!h]
\begin{center}
\includegraphics[scale=1]{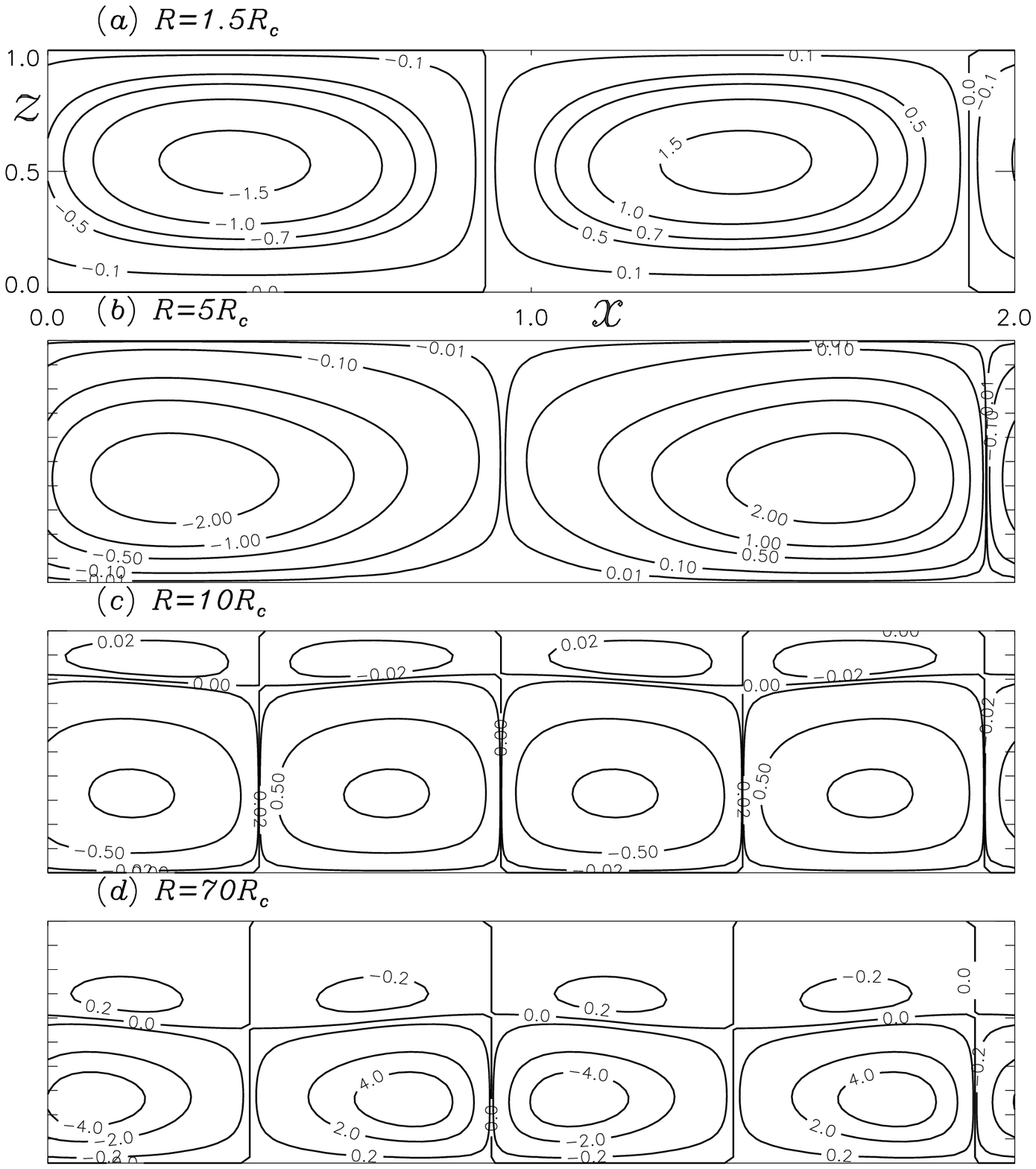}
\end{center}
\caption{Streamlines of steady solutions for $V_c=10$,~$\kappa=1$,~$I_c=0.52$. \label{Fig3}}
\end{figure}

\begin{figure}[!h]
\begin{center}
\includegraphics[scale=1]{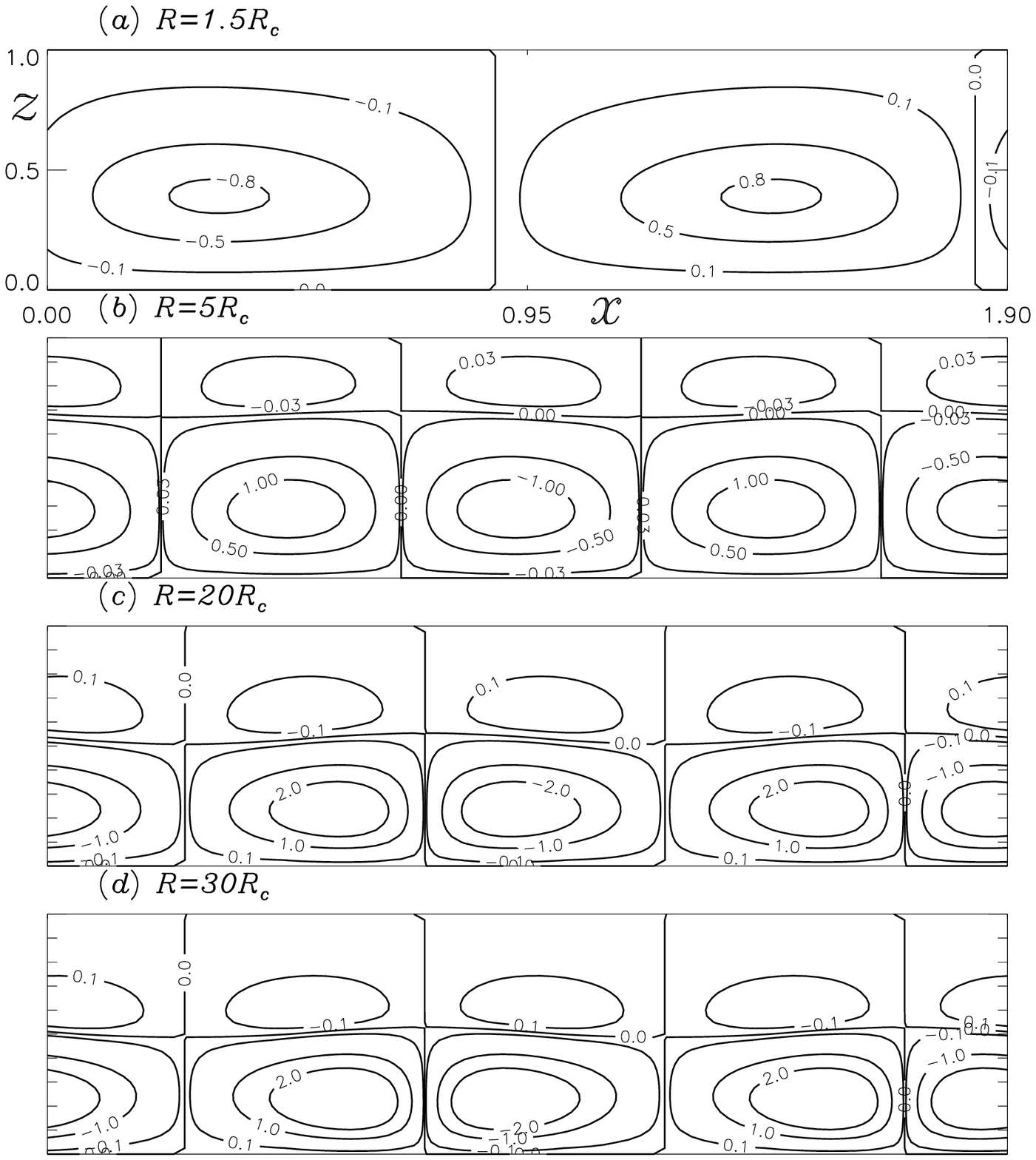}
\end{center}
\caption{Streamlines of steady solutions for $V_c=10$,~$\kappa=1$,~$I_c=0.495$. \label{Fig4}}
\end{figure}

\subsection{Scaling of the equations}
\label{chapt4:scaling}

Length is scaled on the depth $H$, velocity on $D/H$, time on the diffusive
scale ${H^{2}/D}$, and the cell concentration on the mean cell concentration
$\overline{n}$. Now, the governing system for bioconvection in terms of stream function $\psi$ and vorticity $\zeta$, is recasted as:

\begin{equation}
\zeta=-\nabla^{2}\psi,
\label{chapt4:scaled_zeta}
\end{equation}

\begin{equation}
\frac{\partial{\zeta}}{\partial{t}}+{\bmi \nabla}\cdot{\left(\zeta{\bmi u}\right)}=S_{c}\nabla^{2}{\zeta}-S_{c}R\frac{\partial{n}}{\partial{x}}, 
\label{chapt4_scaled_momentum}
\end{equation}

and

\begin{equation}
\frac{\partial{n}}{\partial{t}}=-{\bmi \nabla}\cdot{\bmi J},
\label{chapt4:scaled_cellflux}
\end{equation}

where the flux of the cells is

\begin{equation}
{\bmi J}=n{\bmi u}+nV_{c}<{\bmi p}>-{\bmi \nabla}{n}.
\label{chapt4:nondimensional_flux}
\end{equation}

Here $S_{c}=\nu/D$ is the Schmidt number, $V_{c}=W_{c}H/D$ is the scaled
swimming speed, and $R={\overline{n}v\Delta{\rho{g}H^{3}}}/{\rho{\nu{D}}}$
is the Rayleigh number. Also $<{\bmi p}>=M(I)\hat{\bmi z}$ with $I$ given by

\begin{equation}
I(x,z)=I_{t}\exp\left[-\kappa{\int\limits_{z}^{1}n(x,z)dz}\right], \label{LIS}
\end{equation}

where $\kappa=\alpha{\overline{n}H}$ is the nondimensional extinction (absorption)
coefficient. The phototaxis functions is generated by superimposing the sine functions:

\begin{equation}
M\left(I\right)=0.8\sin{\left[\frac{3\pi}{2}\chi{\left(I\right)}\right]}-0.1\sin{\left[\frac{\pi}{2}\chi{\left(I\right)}\right]}, \label{phototaxis}
\end{equation}

where $\chi{\left(I\right)}=I\exp{\left[\beta\left(I-1\right)\right]}$. Eqs.~(\ref{chapt4:scaled_zeta})--(\ref{chapt4:scaled_cellflux})
are solved in the region $0\leq{x\leq{\lambda}}$~and~$0\leq{z\leq{1}}$, where $\lambda=L/H$~is the normalized width of the domain. 

We impose rigid, no-slip boundary condition on $z=0$ and periodic boundary condition at $x=0, \lambda$ (i.e. $\psi{(0,z)}=\psi{(\lambda,z)}$). We require that both the normal velocity and tangential stress vanish on $z=1$. Also there is no flux of cells through the walls. Thus the boundary conditions are

\begin{eqnarray}
\psi=0 \quad \mbox{and\ }\quad {\bmi J}\cdot{\hat{\bmi z}}=0 \quad \mbox{at\ }\quad z=0,1,\label{scaled_bc1}\\
 \frac{\partial{\psi}}{\partial{z}}=0 \quad \mbox{at}\quad z=0\quad
\mbox{and}\quad \zeta=0\quad \mbox{at}\quad z=1,\label{scaled_bc2}\\
\psi=0,~ \zeta=0 \quad \mbox{and}\quad {\bmi J}\cdot{\hat{\bmi x}=0 \quad \mbox{at}\quad} x=0,\lambda  \label{scaled_bc3}.
\end{eqnarray}

We choose $\lambda=\lambda_c$ to compute the solution in the full convection cell and the initial conditions are $\psi=0,~\zeta=0,~n=1+\epsilon{\cos\left(\pi{x}/\lambda\right)}, \quad \mbox{where} \quad \epsilon=10^{-5}$.

\section{The numerical procedure}
\label{chapt4:numerical_procedure}

The governing Eqs.~(\ref{chapt4:scaled_zeta})--(\ref{chapt4:scaled_cellflux}) with appropriate boundary conditions are solved using a conservative finite-difference scheme via stream function-vorticity formulation \cite{gh,ghorai_thesis} and  the corresponding critical values at the onset of bioconvection are computed in the $(k,R)$--plane by linear stability theory as similar to Panda and Singh \cite{pr}.

\section{Results}
\label{chapt4:results}

The representative parameter values are $V_{c}=10,15, 20$ and $\kappa=0.5,1.0$ respectively. The value of critical intensity $I_{c}$ is selected so that the sublayer where the cells aggregate at equilibrium state lies either around $z=3/4$ or $z=1/2$ of the suspension (see Fig. \ref{basic}). In this investigation, the Rayleigh number is varied such that $1.5R_c \le R \le 100 R_c$ and some steady convection solutions of the proposed study are presented here [see Figs. \ref{Fig1}--\ref{Fig4}].

\section{$V_c=10$}

\subsection{Extinction coefficient $\kappa=0.5$}

We start with the case, when $V_c=10,$ $\kappa=0.5,$ and $I_c=0.66.$ In this case, the sublayer via equilibrium state is located at around three-quarter height of the domain [see Fig. \ref{basic}(a)]. Fig. \ref{Fig1} shows the bioconvective solutions for Rayleigh number $R$ such that $1.5R_c \le R \le 30R_c.$ For $R=1.5R_c,$ a steady state two convection cells solution is observed. When $R$ is increased to $5R_c$ the two convection cells solution is replaced by a four convection cells solution and it persists upto $R=20R_c.$ When $R$ is increased to $30R_c,$ small counter rotating cells appear above the convection cells. If $R$ is increased further, the bioconvective solution becomes periodic.

Next, we consider the case when $V_c=10,$ $\kappa=1,$ and $I_c=0.63.$ In this case, the sublayer via equilibrium state is located at around mid-height of the domain [see Fig. \ref{basic}(a)]. Fig. \ref{Fig2} shows the bioconvective solutions for Rayleigh number $R$ such that $1.5R_c \le R \le 40R_c.$ For $R=1.5R_c,$ a steady state two convection cells solution is observed. When $R$ is increased to $5R_c$ the two convection cells solution is replaced by a four convection cells solution and small counter rotating cells appear above the convection cells. At $R=20R_c,$ the small counter rotating cells which appear above the convection cells grow in height and the trend continues when $R$ is increased to $40R_c.$ If $R$ is increased further, the bioconvective solution becomes periodic.

\subsection{Extinction coefficient $\kappa=1$}

Here start with the case, when $V_c=10,$ $\kappa=1,$ and $I_c=0.52.$ In this case, the sublayer via equilibrium state is located at around three-quarter height of the domain [see Fig. \ref{basic}(b)]. Fig. \ref{Fig3} shows the bioconvective solutions for Rayleigh number $R$ such that $1.5R_c \le R \le 70R_c.$ For $R=1.5R_c,$ a steady state two convection cells solution is observed. When $R$ is increased to $10R_c$ the two convection cells solution is replaced by a four convection cells solution and small counter rotating cells appear above the convection cells. The same trend continues when $R$ is increased upto $70R_c,$ small counter rotating cells which appear above the convection cells grow in height. If $R$ is increased further, the bioconvective solution becomes periodic.

Next, we consider the case when $V_c=10,$ $\kappa=1,$ and $I_c=0.495.$ In this case, the sublayer via equilibrium state is located at around mid-height of the domain [see Fig. \ref{basic}(b)]. Fig. \ref{Fig4} shows the bioconvective solutions for Rayleigh number $R$ such that $1.5R_c \le R \le 30R_c.$ For $R=1.5R_c,$ a steady state two convection cells solution is observed. When $R$ is increased to $5R_c$ the two convection cells solution is replaced by a four convection cells solution and small counter rotating cells appear above the convection cells. At $R=20R_c,$ the small counter rotating cells which appear above the convection cells grow in height and the trend continues when $R$ is increased to $30R_c.$ If $R$ is increased further, the bioconvective solution becomes periodic.

\section{$V_c=15$ and $V_c=20$}
We have also found the bioconvective solutions for $V_c=15$ and $V_c=20$ too and the solutions are qualitatively similar to those of $V_c=10.$

\section{CONCLUSIONS}
\label{chapt4:conclusions}

In this study, two-dimensional phototactic bioconvection in a suspension of non-scattering algae is simulated numerically. Since small domains with periodic sidewalls can be considered as a model of a spatially extended system, the suspension is confined by a stress-free top boundary, and rigid bottom and periodic lateral boundaries. The intensity of light via critical intensity is adjusted so that the sublayer at the equilibrium state lies either at the midheight or three-quarter height of the domain. The non-scattering model for phototaxis proposed by Vincent and Hill~\cite{vh:hv} is employed in this study and the governing bioconvective system is solved using a
conservative finite-difference scheme via stream function-vorticity formulation.

The discrete parameters taken in this study are $V_{c}=10, 15, 20$ and $\kappa=0.5, 1.0$ respectively and the Rayleigh number $R$ is varied such that $1.5R_{c} \le R \le 100R_{c}$. The conclusions drawn from the present study are as follows. The bioconvection solutions always appear in the form of a two-rolls pattern (or any even number of rolls) due to presence of periodic lateral boundaries. The number of convection cells increases as the Rayleigh number increases. The symmetry with respect to the midvertical line in bioconvection cells slution does not persist. Also, the weak counter rotating cells appear on  top of the main convection cells and grow in size as the Rayleigh number increases. 

To test this proposed theoretical model with the available experimental results, we refer the quantitative study by Williams and Bees~\cite{williams_11} and it is in agreement with the experimental results via appearance of multitude of plumes with the increment in Rayleigh number. Penetrative phototactic bioconvection in a non-scattering suspension is in the process of extension to three dimensions since bioconvection is intrinsically a three-dimensional phenomena.
  
\section*{Acknowledgements}

The corresponding author gratefully acknowledges the COUNCIL OF SCIENTIFIC AND INDUSTRIAL RESEARCH (CSIR), Government of India for the financial support via the Extramural Research Grant (Grant No. 25(0295)/19/EMR-II).

\bibliographystyle{plain}

\end{document}